\documentclass[11pt]{amsart}

 \usepackage{amsfonts,graphics,amsmath,amsthm,amsfonts,amscd,amssymb,amsmath,latexsym,multicol,
 mathrsfs}
\usepackage{epsfig,url}
\usepackage{flafter}
\usepackage{fancyhdr}
\usepackage{hyperref}
\hypersetup{colorlinks=true, linkcolor=black}

 \usepackage[matrix, arrow]{xy}

\DeclareMathOperator{\lct}{lct}

\DeclareMathOperator{\mld}{mld}

\DeclareMathOperator{\Supp}{Supp}

\DeclareMathOperator{\vol}{vol}
\DeclareMathOperator{\Bir}{Bir}

 \numberwithin{equation}{subsection}
 \numberwithin{footnote}{subsection}

 \newtheorem{cor}[subsection]{Corollary}

 \newtheorem{thm}[subsection]{Theorem}
 \newtheorem{conj}[subsection]{Conjecture}

{
\theoremstyle{upright}

}

 \newcommand{\N}{\mathbb N}
 \newcommand{\PP}{\mathbb P}
 
 \newcommand{\Q}{\mathbb Q}
 \newcommand{\R}{\mathbb R}
 \newcommand{\Z}{\mathbb Z}
  
 \newcommand{\bir}{\dashrightarrow}
 \newcommand{\rddown}[1]{\left\lfloor{#1}\right\rfloor} 

\title{Birational geometry of algebraic varieties}

\author{Caucher Birkar}
\date{\today}
\begin{document}
\maketitle



\tableofcontents


\section{\bf Introduction}

This is a report on some of the main developments in birational geometry in recent years 
focusing on the minimal model program, Fano varieties, singularities and related topics, in characteristic zero. 
 This is not a comprehensive survey of all advances in birational geometry,
e.g. we will not touch upon the positive characteristic case which is a very active 
area of research.  We will work over an algebraically closed field $k$ of characteristic zero. 
Varieties are all quasi-projective.

Birational geometry, with the so-called minimal model program at its core, aims to classify algebraic varieties
up to birational isomorphism by identifying ``nice" elements in each birational class and then classifying such
elements, e.g study their moduli spaces. Two varieties are birational if they contain isomorphic open subsets. 
In dimension one, a nice element in a birational class is simply a smooth and projective element. 
In higher dimension though there are infinitely many such elements in each class, so picking a 
representative is a very challenging problem. Before going any further lets introduce the canonical 
divisor.

\subsection{Canonical divisor}
To understand a variety $X$ one studies subvarieties and sheaves on it. Subvarieties of codimension one 
and their linear combinations, that is, divisors play a crucial role. Of particular importance is the 
canonical divisor $K_X$. When $X$ is smooth this is the divisor (class) whose associated sheaf 
$\mathcal{O}_X(K_X)$ is the canonical sheaf $\omega_X:=\det \Omega_X$ where $\Omega_X$ is the sheaf 
of regular differential forms. When $X$ is only normal, 
$K_X$ is the closure of the canonical divisor of the smooth locus.
In general, the canonical divisor is the only special non-trivial divisor attached to $X$. 
It plays an important role in algebraic geometry, e.g. in duality theory and Riemann-Roch formula, 
and also in differential and arithmetic geometry. It is a central object in birational geometry.

\vspace{0.2cm}
\emph{Example.}
Assume $X=\PP^d$. Then $K_X\sim -(d+1)H$ where $H\subset \PP^d$ is a hyperplane.

\vspace{0.2cm}
\emph{Example.}
Assume $X\subset \PP^d$ is a smooth hypersurface of degree $r$. Then we have $K_X\sim (-d-1+r)H|_X$ 
where $H\subset \PP^d$ is a hyperplane not containing $X$.

\vspace{0.2cm}
\emph{Example.}
If $X$ is a toric variety, then $K_X\sim -\Lambda$ where $\Lambda$ is the sum of the torus-invariant 
divisors. 

\subsection{Varieties with special canonical divisor}

Let $X$ be a projective variety with ``good" singularities (by this we mean klt or lc singularities 
defined below, see \ref{ss-singularities}). 
\begin{displaymath}
\mbox{We say $X$ is} 
\left\{ \begin{array}{l}
\textrm{\it{Fano} ~~ \hspace{4cm} if $K_X$ is anti-ample}\\
\textrm{\it{Calabi-Yau} ~~ \hspace{3cm} if $K_X$ is numerically trivial}\\
\textrm{\it{canonically polarised} ~~ \hspace{1.4cm}  if $K_X$ is ample}
\end{array} \right. \hspace{5cm}
\end{displaymath} 
Note that here we consider Calabi-Yau varieties in a weak sense, that is, we do not 
require the vanishing $h^i(X,\mathcal{O}_X)=0$ for $0<i<\dim X$ which is usually assumed in other contexts. 
For example, abelian varieties are Calabi-Yau by our definition.

The special varieties just defined are of great importance in algebraic geometry 
(e.g. birational geometry, moduli theory, derived categories), 
differential geometry (e.g. K\"{a}hler-Einstien metrics, stability), arithmetic geometry 
(e.g. existence and density of rational points), 
and mathematical physics (e.g. string theory and mirror symmetry). They 
behave much better  compared to a randomly chosen variety. 

\vspace{0.2cm}
\emph{Example.}
Assume $X$ is a smooth projective curve of genus $g$.  If $g=0$, then $X\simeq \PP^1$ which is Fano. 
If $g=1$, then $X$ is an elliptic curve, hence a Calabi-Yau. If $g\ge 2$, then $X$ is canonically polarised.

\vspace{0.2cm}
\emph{Example.}
Assume $X\subset \PP^d$ is a smooth hypersurface of degree $r$. If $r\le d$, then $X$ is Fano. 
If $r=d+1$, then $X$ is Calabi-Yau. If $r>d+1$, then $X$ is canonically polarised.

\subsection{Minimal model program}
Now we give a brief description of the minimal model program (MMP). Pick a variety $W$.
Using resolution of singularities we can modify $W$ so that it is smooth and projective. 
However, being smooth and projective is not very special as in dimension at least two 
these properties are shared by infinitely many other varieties in the same birational class. 
It is then natural to look for a more special representative.
One of the main aims of birational geometry is to show that we can dismantle $W$
birationally and reconstruct it using canonically polarised, Calabi-Yau, and Fano
varieties. To be more precise we want to establish the following conjecture formulated in its simplest form.

\begin{conj}[Minimal model and abundance]\label{conj-intro-min-model-abund}
Each variety $W$ is birational to a projective variety $Y$ with ``good" singularities such that either 
 \begin{itemize}
\item $Y$ is canonically polarised, or 
 
\item  $Y$ admits a Fano fibration, or 

\item $Y$ admits a Calabi-Yau fibration.  
 \end{itemize}
\end{conj}

In particular, even if $W$ is smooth, $Y$ may be singular. In fact singularity theory is 
an indispensable part of modern birational geometry.

 As the name suggests the
conjecture actually consists of two parts, the \emph{minimal model conjecture} and the \emph{abundance 
conjecture}. The minimal model conjecture essentially says that we can find $Y$ such that $K_Y$ is \emph{nef} meaning 
$K_Y$ intersects every curve non-negatively, or else there is a $K_Y$-negative fibration 
$Y\to Z$ which means we have a \emph{Fano fibration}. The abundance conjecture essentially says that if $Y$ is not canonically 
polarised and if it does not admit a Fano fibration, then it admits a $K_Y$-trivial fibration $Y\to Z$ 
which means we have a \emph{Calabi-Yau fibration}. The minimal model conjecture holds in 
dimension $\le 4$ 
[\ref{Mori-flip}][\ref{Shokurov-log-flips}][\ref{Kawamata-termination}][\ref{Shokurov-pl-flips}][\ref{Shokurov-ordered}] in full generality, and in any dimension 
for varieties of general type [\ref{BCHM}] while the abundance conjecture is
proved in dimension $\le 3$ [\ref{Miyaoka}][\ref{Kaw-abundance}], and in any dimension for varieties of 
general type [\ref{Shokurov-nonvanishing}][\ref{kawamata-bpf}] (also see [\ref{B-lc-flips}] and references therein for 
more results). We should also mention that the \emph{non-vanishing conjecture} which is a special case 
of (a suitable reformulation of) the abundance conjecture 
implies the minimal model conjecture [\ref{B-mmodel}][\ref{B-mmodel-II}]. 

Given a smooth projective $W$, how can we get to $Y$? This is achieved via \emph{running the MMP} which 
is a step by step program making the canonical divisor $K_W$ more positive 
by successively removing or replacing curves along which $K_W$ is not positive. It gives a (conjecturally finite) sequence 
 of birational transformations
$$
  W=W_1  \bir  W_2 \bir  \cdots  \bir W_n=Y
$$ 
consisting of \emph{divisorial contractions}, \emph{flips}, and a last step canonically trivial contraction. 
The required contractions [\ref{Shokurov-nonvanishing}][\ref{kawamata-bpf}] and 
flips [\ref{BCHM}][\ref{HM-flip}] exist. An important ingredient is the \emph{finite generation} of the $k$-algebra 
$$
R=\bigoplus_{m\ge 0} H^0(W,mK_W)
$$ 
in its various forms [\ref{BCHM}][\ref{HM-flip}][\ref{Shokurov-pl-flips}].

A serious issue with the MMP is that we do not know whether it actually stops at some step $W_n$. 
What is not clear is if the MMP 
can produce an infinite sequence of flips. In other words, the minimal model conjecture is reduced to the following. 
 
\begin{conj}[Termination]\label{conj-intro-termination}
There is no infinite sequence of flips.
\end{conj}

The two-dimensional case of the MMP is classical developed in the early 
20th century by Castelnuovo, Enriques, etc. The three-dimensional case (in characteristic zero)
was developed in the 70's-90's through work of many people notably Iitaka, Iskovskikh, Kawamata, Koll\'ar, Mori,
Reid, Shokurov, Ueno, etc. 
The higher dimensional case is still conjectural but a large portion of 
it has been established since the turn of the century by many people including Birkar, Cascini, Hacon, 
M$\rm ^c$Kernan, Shokurov, Xu, etc, involving many difficult problems of 
local and global nature.

\subsection{Pluricanonical systems, Kodaira dimension and Iitaka fibration}

Let $W$ be a smooth projective variety. The space of sections $H^0(W,mK_W)$, for $m\in\Z$, 
and their associated linear systems $|mK_W|$  are of great importance.
When $W$ is one-dimensional the linear system $|K_W|$  determines its geometry to a large extent. 
Indeed the genus $g$ of $W$ is just $h^0(W,K_W)$ which is encoded in $|K_W|$. Moreover, if $g\ge 2$, then 
$|K_W|$ is base point free, and if in addition $W$ is not hyperelliptic, then $|K_W|$ defines 
an embedding of $X$ into a projective space of dimension $g-1$. 
In higher dimension, however,  $|K_W|$ often says little about $W$. One instead needs 
to study $|mK_W|$ for all $m\in \Z$ in order to investigate the geometry of $W$. This leads to the notion of 
\emph{Kodaira  dimension} $\kappa(W)$, an important birational invariant of $W$. This is defined 
to be the maximum of the dimension of the images of $W$ under the maps defined by the linear systems $|mK_W|$ 
for $m>0$.
It takes values in $\{-\infty, 0,1,\dots, \dim X\}$ where the case $-\infty$ corresponds to the 
situation  when $h^0(W,mK_W)= 0$ for every $m>0$. 

Assume $\kappa(W)\ge 0$, that is, $h^0(W,mK_W)\neq 0$ for some $m>0$. 
When $m>0$ is sufficiently divisible, $|mK_W|$ defines a rational fibration $W\bir X$ which is called 
the \emph{Iitaka fibration} of $W$. This is usually defined up to birational equivalence. The dimension of $X$ 
is simply the Kodaira dimension $\kappa(W)$. It is often possible to translate questions about $W$ to 
corresponding questions about $X$. An old problem is the following: 

\begin{conj}\label{conj-iitaka-plurican-systems}
Assume  $\kappa(W)\ge 0$. Then there exists $m\in\N$ depending only on 
$\dim W$ such that $|mK_W|$ defines the Iitaka fibration.
\end{conj}

If $W$ is of \emph{general type}, i.e. if $\kappa(W)=\dim W$, then the conjecture is already 
known  [\ref{HM-nir-bnd}][\ref{Takayama}] (also 
see [\ref{HMX}][\ref{HMX2}] for more recent and more general results). In this case we can take $m$ 
 such that $|mK_{W}|$ defines a birational embedding of $W$ into some projective space. Note that 
$W$ is birational to its canonical model $X$ [\ref{BCHM}] which is a canonically polarised variety and 
understanding $|mK_W|$ is the same as understanding $|mK_X|$.

Now assume  $0\le \kappa(W)<\dim W$. The most general known result is that the conjecture 
is true if we have bounds on certain invariants of the general fibres of the Iitaka fibration [\ref{BZh}]. 
This is done by using a canonical bundle formula for the Iitaka fibration and translating the conjecture 
into a question on the base of the fibration. Very roughly [\ref{BZh}] says that the conjecture 
holds if one understands the case $\kappa(W)=0$. Note that in this case, assuming the minimal model 
and abundance conjectures,  $W$ is birational to a 
Calabi-Yau variety, and understanding $|mK_W|$ is the same as understanding such systems on the 
Calabi-Yau variety.

Finally, assume $\kappa(W)=\infty$. Then all the linear systems $|mK_W|$, for $m>0$, are 
empty.  By the minimal model and 
abundance conjectures, $W$ is birational to a variety $Y$ admitting a Fano fibration $Y\to Z$. 
The general fibres of this fibration are  Fano varieties. It is then 
natural to focus on Fano varieties $F$ and study the linear systems $|-mK_F|$, for $m>0$, in detail. 
There has been extensive studies of these systems, 
especially in low dimension, but general higher dimensional results are quite recent [\ref{B-Fano}][\ref{B-BAB}].

\subsection{Fano varieties, and connection with families, singularities, and termination}

Let $X$ be a Fano variety. A difficulty with investigating 
$|-mK_X|$ is that, unlike the case of varieties of general type, these systems can change dramatically 
if we change $X$ birationally. On the other hand, a standard inductive technique to study $|-mK_X|$ is to use the 
elements of $|-mK_X|$ (usually with bad singularities) to create a particular kind of covering family of  
subvarieties of $X$ and then use induction by restricting to members of this family.  
A  difficulty  in this approach is that  a member of this family is not necessarily Fano, so it is hard to apply induction, 
again unlike the case of varieties of general type.  Despite these difficulties there has been lots of progress 
in recent years.
 
In general there is $m\in\N$ depending only on $\dim X$ such that $|-mK_X|$ is non-empty. 
Moreover, there is an element of $|-mK_X|$ with good singularities [\ref{B-Fano}, Theorem 1.1]: 
this is a special case of \emph{boundedness of complements} (see \ref{t-bnd-compl-usual-local} below). 
In addition if we put a bound on the singularities of $X$, that is, 
if $X$ is $\epsilon$-lc where $\epsilon>0$, then we can choose $m$ so that 
$|-mK_X|$ defines a birational embedding of $X$ into some 
projective space [\ref{B-Fano}, Theorem 1.2] (see \ref{t-eff-bir-e-lc} below). In fact one can go further 
in this case and show that we can choose $m$ so that $-mK_X$ is very ample, hence 
$|-mK_X|$ defines an embedding of $X$ into some projective space, and that the set of such $X$ form a 
bounded family [\ref{B-BAB}, Theorem 1.1]: this is the so-called \emph{BAB conjecture}
(see \ref{t-BAB} below). These results are proved along with various other results and in conjunction 
with Shokurov's theory of complements. We will give ample explanations in subsequent sections. 

So far we have only mentioned \emph{global Fano} varieties but there are other (relative) Fano varieties.  
Assume $X$ has good singularities, $f\colon X\to Z$ is a surjective projective morphism, and $-K_X$ is ample over $Z$. 
We call $X$ \emph{Fano over} $Z$. If $Z$ is a point, then $X$ is a usual Fano variety otherwise in general 
$X$ is not projective. When $\dim X>\dim Z>0$, then 
$f$ is a Fano fibration. Such fibrations appear naturally in birational geometry,
and in other contexts, e.g. families and moduli of Fano's. 

Now assume $f$ is birational. A special 
case is a \emph{flipping contraction}, one of the corner stones of the MMP. Existence of flips basically means 
understanding the  linear systems $|-mK_X|$ relatively over $Z$. Another important special case 
is when $f$ is the identity morphism in which case we are just looking at the {germ of a point} on a 
variety, hence we are doing \emph{singularity theory}. 
Another connection with {singularity theory} is that of singularities of 
\emph{$\R$-linear systems} of divisors on varieties, in general, that is the variety may not be Fano and 
the divisors may not be related to canonical divisors (see \ref{t-bnd-lct} below). This is necessary 
for the proof of BAB. Therefore, studying Fano varieties in the relative setting naturally overlaps with other 
important topics in birational and algebraic geometry.

There is also connection with the termination conjecture. It is understood that the termination 
conjecture is about understanding 
singularities (see \ref{ss-mld-termination}). Moreover, understanding singularities is essentially about understanding 
Fano varieties in the relative birational case. On the other hand, problems about families of Fano 
varieties fits well in this theory (see \ref{ss-McSh-conj} and \ref{ss-stable-Fano}). 
It is then no surprise that recent advances on Fano varieties  
described above is expected to have a profound impact on further developments in birational geometry.

\section{\bf Preliminaries}

In this section we recall some basic notions. We will try to keep technicalities to a minimum 
throughout the text. Most of what we need can be found in [\ref{kollar-mori}][\ref{BCHM}].

\subsection{Contractions}

A \emph{contraction} is a projective morphism $f\colon X\to Z$ of varieties such that 
$f_*\mathcal{O}_X=\mathcal{O}_Z$. In particular, $f$ is surjective with connected fibres.
 
 \subsection{Hyperstandard sets}\label{ss-dcc-sets}

Let $\mathfrak{R}$ be a subset of $[0,1]$. We define 
$$
\Phi(\mathfrak{R})=\left\{1-\frac{r}{m} \mid r\in \mathfrak{R}, m\in \N\right\}
$$
to be the set of \emph{hyperstandard multiplicities} associated to $\mathfrak{R}$. We usually assume  
$0,1\in \mathfrak{R}$ without mention, so $\Phi(\mathfrak{R})$ includes $\Phi({\{0,1\}})$.

\subsection{Divisors and resolutions}
In algebraic geometry Weil divisors usually have integer coefficients. However, in birational geometry 
it is standard practice to consider $\R$-divisors. An $\R$-divisor on a normal variety $X$ is of the form $M=\sum a_iM_i$ where 
$M_i$ are distinct prime divisors and $a_i\in\R$. By $\mu_{M_i}M$ we mean the coefficient $a_i$. 
We say $M$ is $\R$-Cartier if $M$ can be written 
as an $\R$-linear combination of (not necessarily prime) Cartier divisors. For two $\R$-divisors $M$ and $N$, 
$M\sim_\R N$ means $M-N$ is an $\R$-linear combination of principal Cartier divisors (a principal divisor 
is the divisor of zeros and poles of a rational function).

If $X$ is equipped with a 
projective morphism $f\colon X\to Z$, an $\R$-Cartier divisor $M$ is \emph{nef over $Z$} if $M\cdot C\ge 0$ for every curve 
$C$ contracted to a point by $f$. We say $M$ is \emph{ample over $Z$} if it is a positive $\R$-linear combination 
of ample Cartier divisors.
We say $M$ is \emph{big over $Z$} if $M\sim_\R A+D$ where $A$ is ample 
over $Z$ and $D\ge 0$.

A \emph{log resolution} $\phi\colon W\to X$  of $(X,M)$ is a projective birational morphism where 
$W$ is smooth, and the union of the excpetional locus of $\phi$ and the birational transform of $\Supp M$ 
has simple normal crossing singularities.

\subsection{Pairs}\label{ss-singularities}
An important feature of modern birational geometry is that the main objects are pairs rather than 
varieties. Pairs are much better behaved when it comes to induction and passing from a variety 
to a birational model.  

A \emph{pair} $(X,B)$ consists of a normal  variety $X$ and an $\R$-divisor 
$B\ge 0$ such that $K_X+B$ is $\R$-Cartier. 
If the coefficients of $B$ are $\le 1$, we say $B$ is a \emph{boundary}. 

Let $\phi\colon W\to X$ be a log resolution of $(X,B)$. Let 
$$
K_W+B_W:=\phi^*(K_X+B).
$$ 
The \emph{log discrepancy} of a prime divisor $D$ on $W$ with respect to $(X,B)$ 
is defines as 
$$
a(D,X,B):=1-\mu_DB_W.
$$
We say $(X,B)$ is \emph{lc} (resp. \emph{klt})(resp. \emph{$\epsilon$-lc}) 
if every coefficient of $B_W$ is $\le 1$ (resp. $<1$)(resp. $\le 1-\epsilon$). When $B=0$ we just say $X$ is 
lc, etc, instead of $(X,0)$.

A \emph{non-klt place} of $(X,B)$ is a prime divisor $D$ on 
birational models of $X$ such that $a(D,X,B)\le 0$. A \emph{non-klt centre} is the image on 
$X$ of a non-klt place. When $(X,B)$ is lc, a non-klt centre is also called a 
\emph{lc centre}.

If we remove the condition $B\ge 0$, the above definitions still make sense but we add 
\emph{sub} to each notion defined, e.g. instead of lc we say sub-lc, etc.

\vspace{0.2cm}
\emph{Example.} 
The simplest kind of pair is a \emph{log smooth} one, that is, a pair $(X,B)$ where $X$ is smooth 
and $\Supp B$ has simple normal crossing singularities. In this case $(X,B)$ being lc (resp. klt) means every 
coefficient of $B$ is $\le 1$ (resp. $<1$).

\vspace{0.2cm}
\emph{Example.} 
Let $X$ be the cone over a rational curve of degree $n$ (for a more precise definition see the example following 
Theorem \ref{t-BAB}). Then $X$ is klt. But if $X$ is the cone over an elliptic curve, then $X$ is lc but not klt.

\vspace{0.2cm}
\emph{Example.} 
Let $X$ be a klt surface. Let $\phi\colon W\to X$ be the minimal resolution. The exceptional 
curves are all smooth rational curves and they intersect in a special way. There is a 
whose classification of the possible configurations (cf. [\ref{kollar-mori}, Section 4]). Once we know the configuration and 
the self-intersections of the exceptional divisors it is a matter of an easy calculation to 
determine all the log discrepancies.

\subsection{Generalised pairs}\label{ss-gpp}
These pairs appear mainly when one considers the canonical bundle formula 
of a fibration, e.g. see case (2) of \ref{ss-proof-bnd-compl}. A generalised  pair 
is roughly speaking a pair together with a birational polarisation, that is, a nef divisor on some 
birational model. They play an important role in relation with Conjecture \ref{conj-iitaka-plurican-systems} 
[\ref{BZh}] and most of the results of [\ref{B-Fano}][\ref{B-BAB}].
For the sake of simplicity we will try to avoid using these pairs and their subtle properties 
as much as possible but for  convenience here we recall the definition in the projective case only. 
For detailed studies of generalised pairs see [\ref{BZh}][\ref{B-Fano}]. 

A projective \emph{generalised (polarised) pair} consists of 
\begin{itemize}
\item a normal projective variety $X'$,

\item an $\R$-divisor $B'\ge 0$ on $X'$, 

\item a projective birational morphism $\phi\colon X\to X'$ from a normal variety, and 

\item a nef $\R$-Cartier divisor $M$ on $X$, 
\end{itemize}
such that  $K_{X'}+B'+M'$ is $\R$-Cartier, where $M' := \phi_*M$. 
We usually refer to the pair by saying $(X',B'+M')$ is a projective generalised pair with 
data $X\overset{\phi}\to X'$ and $M$. However, we want $\phi$ and $M$ to be birational 
data, that is, if we replace $X$ with a higher model, e.g. a resolution, and replace $M$ with 
its pullback, then we assume the new data defines the same generalised pair. 

Now we define generalised singularities.
Replacing $X$ we can assume $\phi$ is a log resolution of $(X',B')$. We can write 
$$
K_X+B+M=\phi^*(K_{X'}+B'+M')
$$
for some uniquely determined $B$.  We say $(X',B'+M')$ is 
\emph{generalised lc} (resp. \emph{generalised klt})
if every coefficient of $B$ is $\le 1$ (resp. $<1$).

\vspace{0.2cm}
\emph{Example.} 
Assume $(X',B'+M')$ is a  projective generalised pair with 
data $X\overset{\phi}\to X'$ and $M$, and assume $M=\phi^*M'$. 
Then $(X',B'+M')$ is generalised lc (resp. generalised klt) iff $(X',B')$ is lc (resp. klt). 
In other words, in this case $M'$ does not contribute to singularities.

\vspace{0.2cm}
\emph{Example.} 
Let $X'=\PP^2$, and $\phi\colon X\to X'$ be 
the blowup of a closed point $x'\in X'$. Assume $H\subset $ is a hyperplane. 
Let $M=3\phi^*H-tE$ where $E$ is the exceptional divisor of $\phi$ 
and $t\in [0,3]$ is a real number. 
Then letting $B'=0$, $(X',B'+M')$ is a projective generalised pair with data  $X\overset{\phi}\to X'$ and $M$. 
Note that $M'=0$. 
Now we can determine $B$ in the formula above. Calculating intersection numbers 
we find $B=(t-1)E$. Therefore, $(X',B'+M')$ is generalised lc (resp. generalised klt) iff $t\le 2$ (resp. $t<2$).

\section{\bf Fano varieties} 

\subsection{Facets of Fano varieties}

Grothendieck insisted on studying varieties (and schemes) in a relative setting. This philosophy 
has been very successfully implemented in birational geometry. This is particularly interesting 
in the case of Fano varieties, or we should say relative Fano varieties.
 
Let $(X,B)$ be a klt pair and $f\colon X\to Z$ be a surjective projective morphism, and assume
$-(K_X+B)$ is ample over $Z$. We then say $(X,B)$ is \emph{Fano over $Z$}. This relative notion 
unifies various classes of objects of central importance. There are three 
distinct cases. 
\begin{itemize}
\item \emph{Global case:} this is when $Z$ is just a point, hence $(X,B)$ is a \emph{Fano pair} in the usual sense. 

\item \emph{Fibration case:} this is when $\dim X>\dim Z>0$, that is, $f$ is a genuine fibration and its general fibres 
are global Fano pairs. 

\item \emph{Birational case:} this is when $f$ is birational. There are several important subcases here. 
If $f$ is extremal and contracts one divisor,  then $f$ is a \emph{divisorial contraction}. If 
$f$ is extremal and contracts some subvariety but not a divisor, then $f$ is a \emph{flipping contraction}. 
If $f$ is an isomorphism, then $(X,B)$ is just the germ of a klt singularity.

\end{itemize}

\subsection{Complements and anti-pluri-canonical systems}

Assume $(X,B)$ is an lc pair equipped with a projective morphism $X\to Z$. The \emph{theory of complements} 
is essentially the study of the systems $|-n(K_X+B)|$ where $n\in\N$, in a relative sense over $Z$. Obviously 
this is interesting only when some of these systems are not non-empty, e.g. Fano case. The 
theory was introduced by Shokurov [\ref{Shokurov-log-flips}]. The theory 
was further developed in [\ref{shokurov-surf-comp}][\ref{PSh-I}][\ref{PSh-II}][\ref{B-Fano}][\ref{B-BAB}].

A \emph{strong $n$-complement} of $K_{X}+B$ over a point $z\in Z$ is of the form 
$K_{X}+{B}^+$ where over some neighbourhood of $z$ we have: 
\begin{itemize}
\item $(X,{B}^+)$ is lc, 

\item $n(K_{X}+{B}^+)\sim 0$, and 

\item ${B}^+\ge B$.
\end{itemize}
From the definition we get 
$$
-n(K_{X}+B)\sim n{B}^+-nB\ge 0
$$
over some neighbourhood of $z$ which in particular means the linear system $|-n(K_X+B)|$
is not empty over $z$, and that it contains a ``nice" element. 
An $n$-complement [\ref{B-Fano}] is defined similarly but it is more complicated, so for simplicity we avoid using it. 
However, if $B=0$, a complement and a strong complement are the same thing.

\begin{thm}[{[\ref{B-Fano}, Theorems 1.7, 1.8, 1.9]}]\label{t-bnd-compl-usual-local}
Let $d$ be a natural number and $\mathfrak{R}\subset [0,1]$ be a finite set of rational numbers.
Then there exists a natural number $n$
depending only on $d$ and $\mathfrak{R}$ satisfying the following.  
Assume $(X,B)$ is a pair and $X\to Z$ a contraction such that 
\begin{itemize}
\item $(X,B)$ is lc of dimension $d$,

\item the coefficients of $B$ are in $\Phi(\mathfrak{R})$, 

\item $X$ is Fano type over $Z$, and 

\item $-(K_{X}+B)$ is nef over $Z$.
\end{itemize}
Then for any point $z\in Z$, there is a strong $n$-complement $K_{X}+{B}^+$ of $K_{X}+{B}$ 
over $z$. Moreover, the complement is also an $mn$-complement for any $m\in \N$. 
\end{thm}

Here $X$ of \emph{Fano type} over $Z$ means $(X,G)$ is Fano over $Z$ for some $G$.
The theorem was conjectured by Shokurov [\ref{shokurov-surf-comp}, Conjecture 1.3] who proved it in 
dimension $2$ [\ref{shokurov-surf-comp}, Theorem 1.4] (see also [\ref{PSh-II}, Corollary 1.8], 
and [\ref{Shokurov-log-flips}] for some cases). 
Prokhorov and Shokurov [\ref{PSh-I}][\ref{PSh-II}] prove various inductive statements 
regarding complements including some unconditional cases in dimension $3$. 

\vspace{0.2cm}
\emph{Example.}
When $X\to Z$ is toric morphism and $B=0$ we can take $n=1$ and $B^+$ to be the sum of the 
torus-invariant divisors. 

\vspace{0.2cm}
\emph{Remark.} 
Assume $Z$ is a point. 
Assume for simplicity that $B=0$ and that $-K_X$ is ample, that is, $X$ is a 
usual Fano variety. When $X$ is a smooth 3-fold, Shokurov [\ref{Shokurov-anticanonical}] 
proved that $|-K_X|$ contains a 
smooth K3 surface. In particular, $K_X$ has a $1$-complement. This is probably where the 
higher dimensional theory of complements originates. 

\vspace{0.2cm}
\emph{Remark.}  
Assume $X\to Z$ is birational. 
Assume again for simplicity that $B=0$ and that $-K_X$ is ample over $Z$.
When $X\to Z$ is a flipping contraction contracting one smooth rational curve only, Mori (cf. [\ref{kollar-mori-flip}, Theorem 1.7]) 
showed that there always exists a $1$-complement over each $z\in Z$ but in the analytic sense, i.e. it exists over an 
analytic neighbourhood of $z$. This is used in Mori's proof of existence of 
3-fold flips [\ref{Mori-flip}]. 

\vspace{0.2cm}
\emph{Remark.}
Assume $X\to Z$ is an isomorphism, so we are looking at the germ of a klt singularity $(X,B)$ 
around a point $x\in X$. For simplicity again assume $B=0$. In general the Cartier index of 
$K_X$ is not bounded even in dimension $2$. The point of complement theory in this case 
is that the $n$-complement $K_X+B^+$ has Cartier index $n$ which is bounded. 

\vspace{0.2cm}
\emph{Remark.} 
When $X$ is a $3$-fold with terminal singularities, $-K_X$ is ample over $Z$, and $B=0$, 
the \emph{general elephant conjecture} 
of Reid asks whether a general element of the linear system $|-K_X|$, relatively over $Z$, 
has canonical singularities. This is true in various cases, e.g. when $X$ is Gorenstein and $Z$ is a point 
[\ref{Reid-can-Fano}], or when $X\to Z$ is identity [\ref{Reid-term-sing}].

\vspace{0.2cm}
\emph{Example.}
Lets look at the particular case of surfaces in the local case. 
Assume $X$ is a surface, $X\to Z$ is the identity, and $B=0$.
If $x\in X$ is smooth, then $K_X$ is a $1$-complement 
of itself, that is, we can take $B^+=0$. In the singular case, from classification of 
the possible singularities one gets [\ref{shokurov-surf-comp}, 5.2.3]:
\begin{displaymath}
\mbox{if $x\in X$ is a type} 
\left\{ \begin{array}{l}
\textrm{$A$ singularity, then $K_X$ has a $1$-complement.}\\
\textrm{$D$ singularity, then $K_X$ has a $2$-complement.}\\
\textrm{$E_6$ singularity, then $K_X$ has a $3$-complement.}\\
\textrm{$E_7$ singularity, then $K_X$ has a $4$-complement.}\\
\textrm{$E_8$ singularity, then $K_X$ has a $6$-complement.}\\
\end{array} \right. \hspace{2cm}
\end{displaymath}

  \subsection{Effective birationality}   

Let $X$ be a Fano variety. Theorem \ref{t-bnd-compl-usual-local} says that $|-mK_X|$ is non-empty 
containing a nice element for some $m>0$ depending only on $\dim X$. 
If we bound the singularities of $X$, we then have a much stronger statement.

\begin{thm}[{[\ref{B-Fano}, Theorem 1.2]}]\label{t-eff-bir-e-lc}
Let $d$ be a natural number and $\epsilon>0$ a real number. Then there is a natural 
number $m$ depending only on $d$ and $\epsilon$ such that if $X$ is any $\epsilon$-lc Fano 
variety of dimension $d$, then $|-mK_X|$ defines a birational map.
\end{thm}

Note that $m$ indeed depends on $d$ as well as $\epsilon$ because the theorem implies 
the volume $\vol(-K_X)$ is bounded from below 
by $\frac{1}{m^d}$. Without the $\epsilon$-lc assumption, $\vol(-K_X)$ can get arbitrarily small 
or large [\ref{HMX2}, Example 2.1.1]. In dimension $2$, the theorem is a consequence of BAB [\ref{Alexeev}], 
and in dimension $3$, special cases are proved in [\ref{Jiang}] using different methods. 
Cascini and M$^{\rm{c}}$Kernan have independently proved the 
theorem for canonical singularities, that is when $\epsilon=1$, using quite different methods.

It is worth mentioning that the theorem also holds in the relative setting. It follows immediately from the 
global case stated above.\\

\subsection{Boundedness of Fano varieties: BAB}

It is possible to strengthen \ref{t-eff-bir-e-lc} so that $|-mK_X|$ defines an actual embedding. 
This follows from the next result.

\begin{thm}[{[\ref{B-BAB}, Theorem 1.1]}]\label{t-BAB}
Let $d$ be a natural number and $\epsilon$ a positive real number. Then the projective 
varieties $X$ such that  
\begin{itemize}
\item $(X,B)$ is $\epsilon$-lc of dimension $d$ for some boundary $B$, and 

\item $-(K_X+B)$ is nef and big,
\end{itemize}
form a bounded family. 
\end{thm}

This was known as the Borisov-Alexeev-Borisov or BAB conjecture. 
Various special cases of it was considered by many people. It was known 
in the following cases (by taking $B=0$): surfaces [\ref{Alexeev}], toric varieties [\ref{A-L-Borisov}], 
Fano $3$-folds with terminal singularities and Picard number one [\ref{kawamata-Fano-term-3-folds}],
 Fano $3$-folds with canonical singularities  [\ref{KMMT-can-3-folds}], smooth Fano varieties [\ref{KMM-smooth-fano}],
 spherical Fano varieties [\ref{Alexeev-Brion}], 
 Fano $3$-folds with fixed 
 Cartier index of $K_X$ [\ref{A-Borisov}], and more generally, Fano varieties of given dimension with fixed
 Cartier index of $K_X$ [\ref{HMX2}]; 
in a given dimension, the Fano varieties $X$ equipped with a boundary $\Delta$ such that 
$K_X+\Delta\equiv 0$, $(X,\Delta)$ is $\epsilon$-lc, and such that the coefficients of 
$\Delta$ belong to a DCC set [\ref{HMX2}] (also see [\ref{HX}][\ref{B-Fano}, Theorem 1.4]).

\vspace{0.2cm}
\emph{Example.}
Now we look at an example of a non-bounded family of singular Fano surfaces. For $n\ge 2$ consider 
$$
\xymatrix{
E & \subset   W_n  \ar[d] \ar[r]^{f} & X_n\\
 & \mathbb{P}^1 & }
$$
where $X_n$ is the cone over a rational curve of deg $n$,  $f$ is blowup of the vertex, and $E$ is the exceptional curve. 
In other words, $W_n$ is the projective bundle of $\mathcal{O}_{\mathbb{P}^1}\oplus \mathcal{O}_{\mathbb{P}^1}(-n)$,
$E$ is the section given by the summand $\mathcal{O}_{\mathbb{P}^1}(-n)$, and
$X_n$ is obtained from $W_n$ by contracting $E$. Then an easy calculation, using $E^2=-n$, shows that 
$$
K_{W_n}+\frac{n-2}{n}E=f^*K_{X_n},
$$ 
hence $X_n$ is a $\frac{2}{n}$-lc Fano  variety with one singular point (the larger is $n$, the deeper is the singularity).
In particular, since the set of numbers $\{\frac{n-2}{n} \mid n\in \mathbb{N} \}$ is not finite, 
the set $\{X_n \mid n\in \mathbb{N}\}$ is not a bounded family. This example explains the 
role of the number $\epsilon$ in Theorem \ref{t-BAB}.

\vspace{0.2cm}
\emph{Example.}
In this example we sketch the proof of Theorem \ref{t-BAB} in dimension two following [\ref{Alexeev-Mori}].
For simplicity assume $B=0$ and that $-K_X$ is ample.
There is $\Delta\ge 0$ such that $(X,\Delta)$ is $\epsilon$-lc and $K_X+\Delta\sim_\R 0$. 
Let $\phi\colon W\to X$ be the minimal resolution and let $K_W+\Delta_W$ be the pullback of $K_X+\Delta$. 
Since $(X,\Delta)$ is klt, the exceptional divisors of $\phi$ are all smooth rational curves. 
Moreover, by basic properties of minimal resolutions, $\Delta_W\ge 0$. In particular, $(W,\Delta_W)$ is an 
$\epsilon$-lc pair. Now a simple calculation of intersection numbers shows that $-E^2\le l$ 
for every exceptional curve of $\phi$ where $l\in\N$ depends only on $\epsilon$. If the number 
of exceptional curves of $\phi$ is bounded, then the Cartier index of $-K_X$ is bounded which in turn 
implies $-nK_X$ is very ample for some bounded $n$. In particular, this holds if 
the Picard number of $W$ is bounded from above. If in addition $\vol(-K_X)$ is bounded, then 
$X$ belongs to a bounded family. Note that $\vol(-K_X)=\vol(-K_W)$. 

Running an MMP on $K_W$ we get a morphism $W\to V$ where $V$ is either $\PP^2$ or a rational ruled surface 
(like $W_n$ in the previous example), 
and the morphism is a sequence of blowups at smooth points. Let $\Delta_V$ be the pushdown of $\Delta_W$. 
Then $(V,\Delta_V)$ is $\epsilon$-lc and $K_V+\Delta_V\sim_\R 0$. It is easy to show that there are 
finitely many possibilities for $V$.
In particular, from $\vol(-K_W)\le \vol(-K_V)$, we deduce that $\vol(-K_X)=\vol(-K_W)$ is bounded from above.
Thus it is enough to prove that the number of blowups in $W\to V$ is bounded. This number can be 
bounded by an elementary analysis of possible intersection numbers in the sequence 
(see [\ref{Alexeev-Mori}, Section 1] for more details).

\subsection{Birational automorphism groups}

An interesting consequence of Theorem \ref{t-BAB} concerns the Jordan property of birational 
automorphism groups of rationally connected varieties. Prokhorov and Shramov [\ref{Prokhorov-Shramov}, Theorem 1.8] 
proved the next result assuming Theorem \ref{t-BAB}.

\begin{cor}[{[\ref{B-BAB}, Corollary 1.3]}]\label{cor-bir-aut}
Let $d$ be a natural number. Then there is a natural number $h$ depending only on $d$ satisfying the following. 
Let $X$ be a rationally connected variety of dimension $d$ over $k$. Then for any finite subgroup $G$ of 
the birational automorphism group $\Bir(X)$, there is a normal abelian subgroup 
$H$ of $G$ of index at most $h$. In particular, $\Bir(X)$ is Jordan.
\end{cor}

Here $X$ \emph{rationally connected} means that every two general closed points can be joined by a  rational 
curve. If we take $X=\PP^d$ in the corollary, then we deduce that the Cremona group 
${\rm Cr}_d:=\Bir(\PP^d)$ is Jordan, answering a question of Serre [\ref{Serre}, 6.1].\\

\section{\bf Singularities of linear systems}

\subsection{Lc thresholds of $\R$-linear systems.}
Let $(X,B)$ be a pair. The \emph{log canonical threshold} (lc threshold for short) of 
an $\R$-Cartier $\R$-divisor $L\ge 0$ with respect to $(X,B)$ 
is defined as 
$$
\lct(X,B,L):=\sup\{t \mid (X,B+tL) ~~\mbox{is lc}\}.
$$ 
It is a way of measuring the singularities of $L$ taking into account the singularities of $(X,B)$ as well.

Now let $A$ be an $\R$-Cartier $\R$-divisor. The $\R$-linear system of $A$ is 
$$
|A|_\R=\{L\ge 0 \mid L\sim_\R A\}.
$$
We then define the \emph{lc threshold} of $|A|_\R$ with respect to $(X,B)$ (also called global 
lc threshold or $\alpha$-invariant)  as 
$$
\lct(X,B,|A|_\R):=\inf\{\lct(X,B,L) \mid L\in |A|_\R\}
$$
which coincides with 
$$
\sup\{t \mid (X,B+tL) ~~\mbox{is lc for every} ~~L\in |A|_\R\}.
$$
This is an asymptotic invariant, so not surprisingly it is hard to compute in specific cases and study 
in general.

Due to connections 
with the notion of stability and existence of K\"{a}hler-Einstein metrics, lc thresholds of $\R$-linear systems  
have attracted a lot of attention, particularly, 
when $A$ is ample. An important special case is when $X$ is Fano and $A=-K_X$ in which case 
many examples have been calculated, e.g. see [\ref{cheltsov-shramov}].

\vspace{0.2cm}
\emph{Example.}
If $X=\PP^d$, $B=0$, and $A=-K_X$, then 
$$
\lct(X,B,|A|_\R)=\frac{1}{d+1}.
$$ 
On the other hand, if 
 $X\subset \PP^d$ is a smooth hypersurface of degree $r\le d$, $B=0$, and $A=-K_X$, then 
$$
\lct(X,B,|A|_\R)=\frac{1}{d+1-r}
$$
[\ref{cheltsov-shramov}, Example 1.3].  

Another reason for studying the above threshold is connection with boundedness of Fano varieties. 
Indeed it plays a central role in the proof of Theorem \ref{t-BAB}.

\begin{thm}[{[\ref{B-BAB}, Theorem 1.4]}]\label{t-bnd-lct-global}
Let $d$ be a natural number and $\epsilon$ a positive real number. Then there is a 
positive real number $t$ depending only on $d,\epsilon$ satisfying the following. 
Assume 
\begin{itemize}
\item $(X,B)$ is a projective $\epsilon$-lc pair of dimension $d$, and

\item $A:=-(K_X+B)$ is nef and big.
\end{itemize}
Then 
$$
\lct(X,B,|A|_\R)\ge t.
$$
\end{thm}

This was conjectured by Ambro [\ref{Ambro}] who proved it in the toric case. It can be 
derived from \ref{t-BAB} but in reality it is proved before \ref{t-BAB} (see next section). 
Jiang [\ref{Jiang}][\ref{Jiang-2}] proved it in dimension two. 

The lc threshold of an $\R$-linear system $|A|_\R$ is defined as an infimum of usual lc thresholds. 
Tian [\ref{Tian}, Question 1] asked whether the infimum is a minimum when $A=-K_X$ and $X$ is Fano.  
The question was reformulated and generalised to Fano pairs in [\ref{cheltsov-shramov}, 
Conjecture 1.12].  The next result gives a positive 
answer when the lc threshold is at most $1$.

\begin{thm}[{[\ref{B-BAB}, Theorem 1.5]}]\label{t-global-lct-attained}
Let $(X,B)$ be a projective klt pair such that $A:=-(K_X+B)$ is nef and big. 
Assume that $\lct(X,B,|A|_\R)\le 1$.  Then there is $0\le D\sim_\R A$ such that 
$$
\lct(X,B,|A|_\R)=\lct(X,B,D).
$$
Moreover, if $B$ is a $\Q$-boundary, then we can choose $D\sim_\Q A$, hence 
in particular, the lc threshold is a rational number.
\end{thm}

Shokurov has an unpublished proof of the theorem in dimension two.

\subsection{Lc thresholds of $\R$-linear systems with bounded degree.}
Next we treat lc thresholds associated with divisors on varieties, in a general setting. 
To obtain any useful result, one needs to impose certain boundedness conditions on the invariants of the 
divisor and the variety. 

\begin{thm}[{[\ref{B-BAB}, Theorem 1.6]}]\label{t-bnd-lct}
Let $d,r$ be natural numbers and $\epsilon$ a positive real number. 
Then  there is a positive real number $t$ depending only on $d,r,\epsilon$ satisfying the following. 
Assume 
\begin{itemize}
\item $(X,B)$ is a projective $\epsilon$-lc pair of dimension $d$, 

\item $A$ is a very ample divisor on $X$ with $A^d\le r$, and

\item $A-B$ is ample.
\end{itemize}
Then  
$$
 \lct(X,B,|A|_\R)\ge t.
$$
\end{thm}

\vspace{0.3cm}
This is one of the main ingredients of the proof of Theorem \ref{t-bnd-lct-global} but it is also interesting on its own. 
We explain briefly some of the assumptions of the theorem. The condition $A^d\le r$ means that 
$X$ belongs to a bounded family of varieties, actually, if we choose $A$ general in its linear system, 
then $(X,A)$ belongs to a bounded family of pairs. We can use the divisor $A$ to measure how ``large" 
other divisors are on $X$. Indeed, 
the ampleness of $A-B$  roughly speaking says 
that the ``degree" of  $B$ is bounded from above, that is,
$$
\deg_AB:=A^{d-1}B< A^d\le r.
$$ 
Without such boundedness assumptions, one would not find a positive lower bound for the lc threshold 
as the next example shows.

\vspace{0.2cm}
\emph{Example.}
Assume $(X=\PP^2,B)$ is $\epsilon$-lc and $S\subset X$ is a line. 
Let $L=A=lS$ where $l\in\N$. Then the multiplicity of $L$ at any closed point $x\in L$ is $l$, hence the 
lc threshold $\lct(L,X,B)\le \frac{1}{l}$. Thus the larger is $l$, the smaller is the threshold. Next we illustrate how 
the threshold depends on the degree of $B$. Let $T$ be another line and $x$ be the intersection point $S\cap T$. 
Let $X_1\to X$ be the blowup at $x$, and let $x_1$ be the intersection of the exceptional 
divisor $E_1$ and the birational transform $S^\sim$. Let $X_2\to X_1$ be the blowup at 
$x_1$, and let $x_2$ be the intersection of the new exceptional divisor $E_2$ and $S^\sim$. 
At each step we blowup the intersection point of $S^\sim$ and the newest 
exceptional divisor. 

Put $W:=X_r$. Then the exceptional locus of $\phi\colon W\to X$ consists of a chain of curves all of which 
are $-2$-curves except one which is a $-1$-curve. 
Then $-K_W$ is nef over $X$, in fact, it is semi-ample over $X$. 
Thus there is $0\le B_W\sim_\R \alpha \phi^*H-K_W$ for some $\alpha>0$  
such that $(W,B_W)$ is $\frac{1}{2}$-lc and $K_W+B_W\sim_\R 0/X$. Now let $B$ be the pushdown 
of $B_W$.  Then $(X,B)$ is $\frac{1}{2}$-lc. Now let $L=S+T$. Then the coefficient of $E_r$ in 
$\phi^* L$ is $r+1$, hence 
$$
\lct(L,X,B)=\lct(\phi^* L,W,B_W)\le \frac{1}{r+1}.
$$ 
Thus there is no lower bound on the lc threshold if $r$ is arbitrarily large. 
This does not contradict Theorem \ref{t-bnd-lct} because when $r\gg0 $, the degree $\deg_AB\gg 0$ and 
$A-B$ cannot be ample (here $A=lS$ with $l$ fixed).

\section{\bf Brief sketch of proofs of main results}\label{s-proof}

In this section we sketch some of the ideas of the proofs of Theorems \ref{t-bnd-compl-usual-local}, \ref{t-eff-bir-e-lc}, 
\ref{t-bnd-lct}, \ref{t-BAB}. We try to remove technicalities as much as possible but this comes at 
the expense of being imprecise in various places and not elaborating on many of the new ideas.

\subsection{Sketch of proof of boundedness of complements}\label{ss-proof-bnd-compl}(Theorem \ref{t-bnd-compl-usual-local})
For simplicity we look at the global case, that is, when $Z$ is a point. 
Pick a sufficiently small $\epsilon\in(0,1)$. 
Let $Y\to X$ be the birational morphism which extracts all the prime divisors with log discrepancy 
smaller than $\epsilon$. Let $K_Y+B_Y$ be the pullback of $K_X+B$. Define $\Theta_Y$ 
to be the same as $B_Y$ except that we replace each coefficient in $(1-\epsilon,1)$ with $1$.  
Run an MMP on $-(K_Y+\Theta_Y)$ and let $Y'$ be the resulting model and $\Theta_{Y'}$ 
be the pushdown of $\Theta_Y$. We can run such MMP 
because $Y$ turns out to be of Fano type, so we can run MMP on any divisor on $Y$. 

As a consequence of local and global ACC [\ref{HMX2}, Theorems 1.1 and 1.5] (in practice we 
need their generalisations to generalised pairs [\ref{BZh}, Theorems 1.5 and 1.6]), we can show that 
the MMP does not contract any component of $\rddown{\Theta_Y}$, $(Y',\Theta_{Y'})$ is lc, 
and $-(K_{Y'}+\Theta_{Y'})$ is nef. It is enough to construct a bounded complement for 
$K_{Y'}+\Theta_{Y'}$. Replacing $(X,B)$ with $(Y',\Theta_{Y'})$ and applying further reductions,  
we can reduce the problem to one of the following cases:

\begin{enumerate}
\item $B$ has a component $S$ with coefficient $1$ and $-(K_X+B)$ is nef and big, or 

\item $K_X+B\equiv 0$ along a fibration $f\colon X\to T$, or 

\item $(X,B)$ is \emph{exceptional}.
\end{enumerate}

Here exceptional means that for any choice of $0\le P\sim_\R -(K_X+B)$ the pair $(X,B+P)$ is klt.  
These cases require very different inductive treatment.

Case (1): First apply \emph{divisorial adjunction} to define $K_S+B_S=(K_X+B)|_S$. 
Further modification of the setting allows us to 
ensure that $S$ is Fano type. Moreover, the coefficients of $B_S$ happen to be in a set $\Phi(\mathfrak{S})$ 
for some fixed finite set $\mathfrak{S}$. By induction on dimension $K_S+B_S$ has a strong $n$-complement for 
some bounded $n$. The idea then is to lift the complement to $X$ using vanishing theorems.
In the simplest case when $(X,B)$ is log smooth and ${B}=S$, we look at the exact sequence 
$$
H^0(-n(K_X+B))\to H^0( -n(K_X+B)|_S)\to H^1(-n(K_X+B)-S)=0
$$
where the vanishing follows from Kawamata-Viehweg vanishing theorem noting that 
$$
-n(K_X+B)-S=K_X-n(K_X+B)-(K_X+B)=K_X-(n+1)(K_X+B)
$$ 
Since $K_S+B_S$ has a strong $n$-complement, the middle space in the above sequence 
is non-trivial which implies the left hand side is also non-trivial by lifting the section corresponding 
to the complement. One then argues that the lifted section gives a 
 strong $n$-complement for $K_X+B$. 

Case (2):  Apply the {\it{canonical bundle formula}} (also called adjunction for fibre spaces, 
derived from [\ref{kaw-subadjuntion}]) to  write 
$$
K_X+B\sim_\R f^*(K_T+B_T+M_T)
$$
where $B_T$ is the \emph{discriminant divisor} and $M_T$ is the \emph{moduli divisor}. It turns out that 
the coefficients of $B_T$ happen to be in a set $\Phi(\mathfrak{S})$ 
for some fixed finite set $\mathfrak{S}$, and that $pM_T$ is integral for some bounded number $p\in \N$. 
Now we want to find a complement for $K_T+B_T+M_T$ and pull it back to $X$. 
There is a serious issue here: $(T,B_T+M_T)$ is not a pair in the usual sense but it is a 
generalised pair. Thus we actually need to prove \ref{t-bnd-compl-usual-local} (at least in the 
global case) in the more general setting of generalised pairs. This makes life a lot more difficult 
but fortunately everything turns out to work. Once we have a bounded complement for 
$K_T+B_T+M_T$ it is straightforward to derive a bounded complement for $K_X+B$.

Case (3):  In this case we use effective birationality. Perhaps after decreasing $\epsilon$, 
the exceptionality condition implies that $(X,B)$ is $\epsilon$-lc.  
For simplicity assume $B=0$ and that $X$ is a Fano variety. Also 
assume we already have Theorem \ref{t-eff-bir-e-lc}. Then  
there is a bounded number $m\in\N$ such that $|-mK_X|$ defines a birational map. Pick 
$M\in |-mK_X|$ and let $B^+=\frac{1}{m}M$. Since $X$ is exceptional, $(X,B^+)$ is automatically klt, 
hence  $K_X+B^+$ is a strong $m$-complement. Although this gives some ideas of how 
one may get a bounded complement but in practice we cannot give a complete 
proof of Theorem \ref{t-eff-bir-e-lc} before proving \ref{t-bnd-compl-usual-local}. The two theorems are 
actually proved together. See [\ref{B-Fano}, Sections 6 and 7] for more details.

\subsection{Sketch of proof of effective birationality}(Theorem \ref{t-eff-bir-e-lc})
Let $m\in\N$ be the smallest number such that $|-mK_X|$ defines a birational map, and 
let $n\in\N$ be a number such that ${\rm vol}(-nK_X)>(2d)^d$. Initially we take $n$ to be the smallest 
such number but  we will modify it during the proof. We want to show that $m$ is bounded from above.
The idea is first to show that $\frac{m}{n}$ is bounded from above, and then at the end show that 
$m$ is bounded. 

Applying a standard elementary technique we can create a covering family $\mathcal{G}$ of subvarieties of $X$ such that 
if $x,y\in X$ are any pair of general closed points, then there is $0\le \Delta \sim_\Q -(n+1)K_X$ 
and $G\in\mathcal{G}$ such that $(X,\Delta)$ is lc at $x$ with the unique non-klt centre $G$, 
and $(X,\Delta)$ is not klt at $y$. 

Assume $\dim G=0$ for all $G$. Then $G=\{x\}$ is an isolated 
non-klt centre. Using multiplier ideals and vanishing theorems we can lift sections from $G$ and show  that 
$|-nK_X|$ defines a  birational map after replacing $n$ with a bounded multiple, hence in particular 
$\frac{m}{n}$ is bounded from above in this case. 

Now lets assume all $G$ 
have positive dimension. If ${\rm vol}(-mK_X|_{G})$ is large, then again using some 
elementary arguments, we can  create a new non-klt centre $G'$ containing $x$ but with $\dim G'<\dim G$. 
Thus we can replace $G$ with $G'$ and apply induction on dimension of $G$. We can then assume 
${\rm vol}(-mK_X|_{G})$ is bounded from above. 

Similar to the previous paragraph,  we can cut $G$ and decrease its dimension if 
${\rm vol}(-nK_X|_{G})$ is bounded from below. Showing this lower boundedness 
is the hard part. A key point here is that although $G$ is not necessarily a divisor and 
although the singularities of $(X,\Delta)$ away from 
$x$ maybe quite bad but still there is a kind of adjunction formula, that is, if $F$ is the normalisation of $G$, 
then we can write  
$$
(K_X+\Delta)|_F\sim_\R K_F+\Theta_F+P_F
$$ 
where $\Theta_F$ is a boundary divisor with coefficients in a fixed DCC set $\Psi$ depending 
only on $d$, and  $P_F$ is pseudo-effective. Replacing $n$ with $2n$ and adding to $\Delta$ 
we can easily make $P_F$ big and effective. 

Now we would ideally want to apply induction on $d$ but the difficulty is that $F$ may not be Fano, in fact, it 
can be any type of variety. Another issue is that the singularities of $(F,\Theta_F+P_F)$ can be pretty bad. To 
overcome these difficulties we use the fact that ${\rm vol}(-mK_X|_{G})$ is bounded from above. 
From this boundedness one can deduce that there is a bounded projective log 
smooth pair $(\overline{F},\Sigma_{\overline{F}})$
 and a birational map $\overline{F}\bir F$ such that $\Sigma_{\overline{F}}$ is reduced 
containing the exceptional divisor of $\overline{F}\bir F$ and the support of the birational 
transform of $\Theta_F$ (and other relevant divisors).

Surprisingly, the worse the 
singularities of $(F,\Theta_F+P_F)$ the better because we can then produce divisors on $\overline{F}$ 
with bounded ``degree" but with arbitrarily small lc thresholds which would contradict a baby version of 
Theorem \ref{t-bnd-lct}.
Indeed assume $(F,\Theta_F+P_F)$ is not klt. A careful study of the above adjunction formula allows to 
write $K_F+\Lambda_F:=K_X|_F$ where $\Lambda_F\le \Theta_F$ and $(F,\Lambda_F)$ is sub-$\epsilon$-lc. 
Put $I_F=\Theta_F+P_F-\Lambda_F$. Then 
$$
I_F=K_F+\Theta_F+P_F-K_F-\Lambda_F\sim_\R (K_X+\Delta)|_F-K_X|_F=\Delta|_F\sim_\R -(n+1)K_X|_F.
$$
Moreover, $K_F+\Lambda_F+I_F$ is ample. 

Let $\phi\colon F'\to F$ and $\psi\colon F'\to \overline{F}$ be a common resolution.  
Pull back $K_F+\Lambda_F+I_F$ to $F'$ and then push it down to $\overline{F}$ and write it as 
$K_{\overline{F}}+\Lambda_{\overline{F}}+I_{\overline{F}}$. Then the above ampleness gives 
$$
\phi^*(K_F+\Lambda_F+I_F)\le \psi^*(K_{\overline{F}}+\Lambda_{\overline{F}}+I_{\overline{F}})
$$
which implies that $(\overline{F},\Lambda_{\overline{F}}+I_{\overline{F}})$ is not sub-klt. From this one 
deduces that $(\overline{F},\Gamma_{\overline{F}}+I_{\overline{F}})$ is not klt where 
$\Gamma_{\overline{F}}=(1-\epsilon)\Sigma_{\overline{F}}$. Finally, one argues that the degree 
of $I_{\overline{F}}$ gets arbitrarily small if ${\rm vol}(-nK_X|_{G})$ gets arbitrarily small, and 
this contradicts an easy case of Theorem \ref{t-bnd-lct}.

If singularities of $(F,\Theta_F+P_F)$ are good, then we again face some serious difficulties. 
Very roughly, in this case, we lift sections from $F$ to $X$ and use this section 
to modify $\Delta$ so that  $(F,\Theta_F+P_F)$ 
has bad singularities, hence we reduce the problem to the above arguments. 
This shows $\frac{m}{n}$ is bounded.

Finally, to we still need to bound $m$. This can be done by arguing that $\vol(-mK_X)$ is bounded 
from above and use this to show $X$ is birationally bounded, and then work on the bounded model.
See [\ref{B-Fano}, Section 4] for more details.

\subsection{Sketch of proof of boundedness of lc thresholds} (Theorem \ref{t-bnd-lct})
Pick $0\le N\sim_\R A$. Let $s$ be the largest number such that $(X,B+sN)$ is 
$\epsilon'$-lc where $\epsilon'=\frac{\epsilon}{2}$. It is enough to show $s$ is bounded from below.
There is a prime divisor $T$ on birational models of $X$  with log discrepancy 
$a(T,X,\Delta)=\epsilon'$ where $\Delta:=B+sN$. It is enough to show that the multiplicity of $T$ in $\phi^*N$ is bounded 
on some resolution $\phi\colon V\to X$ on which $T$ is a divisor. We can assume 
the image of $T$ on $X$ is a closed point $x$ otherwise we can cut by hyperplane sections and 
apply induction on dimension. 

There is a birational morphism $Y\to X$ from a normal projective variety which contracts 
exactly $T$. A key ingredient here is provided by the theory of complements: using the fact that 
$-(K_Y+T)$ is ample over $X$, we can find  $\Lambda_Y$ such that $(Y,\Lambda_Y)$ is 
lc near $T$ and $n(K_Y+\Lambda_Y)\sim 0/X$ for some bounded number $n\in\N$. 
One can think of $K_Y+\Lambda_Y$ as a local-global type of complement. 
The crucial point is that if $\Lambda$ is the pushdown of $\Lambda_Y$, then we can make 
sure degree of $\Lambda$ is bounded from above, that is, after replacing $A$ we can assume 
$A-\Lambda$ is ample. By construction,  the log discrepancy 
$a(T,X,\Lambda)=0$ and $(X,\Supp \Lambda)$ is bounded. 

Next using resolution of singularities we can modify the setting and then assume that 
 $(X,\Lambda)$ is log smooth and $\Lambda$ is reduced. The advantage of having 
 $\Lambda$ is that now $T$ can be obtained by a sequence of blowups which is toroidal 
 with respect to  $(X,\Lambda)$. That is, in every step we blowup the centre of $T$ 
which happens to be a stratum of $(X,\Lambda)$; a stratum is just a 
component of the intersection of some of the components of $\Lambda$.
The first step is just the blowup of $x$. One argues that it is enough to bound the number of these blowups.

By the previous paragraph,  we can discard any component of 
$\Lambda$ not passing through $x$, hence assume $\Lambda=S_1+\dots+S_d$ 
where $S_i$ are irreducible components.
On the other hand, a careful analysis of $Y\to X$ allows us to further modify the situation so that 
$\Supp \Delta$ does not contain any stratum of $(X,\Lambda)$ apart from $x$. This is one of the difficult 
steps of the whole proof.

Since $(X,\Lambda)$ is log smooth and bounded, we can find 
a surjective finite morphism $X\to \mathbb{P}^d$ which maps $x$ to the origin $z=(0:\cdots:0:1)$ 
and maps $S_i$ on $H_i$ where $H_1,\dots,H_d$ are the coordinate hyperplanes passing through 
$z$. Since $\Supp \Delta$ does not contain any stratum of $(X,\Lambda)$ apart from $x$, 
it is not hard to reduce the problem to a similar problem on $\mathbb{P}^d$. From now on we assume 
$X=\mathbb{P}^d$ and that $S_i$ are the coordinate hyperplanes. The point of this reduction is that 
now $(X,\Lambda)$ is not only toroidal but actually toric, and $-(K_X+\Lambda)$ is very ample. 
In particular, replacing $\Delta$ with $t\Delta+(1-t)\Lambda$ for some 
sufficiently small $t>0$ (and replacing $\epsilon'$ accordingly), we can make 
$K_X+\Delta$ anti-ample. Next by adding to $\Delta$ we can assume 
$K_X+\Delta$ is numerically trivial.

Let $W\to X$ be the sequence of blowups which obtains $T$ as above.
Since the blowups are toric, $W$ is a toric variety. If $Y\to X$ is the birational morphism 
contracting $T$ only, as before, then $Y$ is also a toric variety. Moreover, if $K_Y+\Delta_Y$ is the 
pullback of $K_X+\Delta$, then $(Y,\Delta_Y)$ is $\epsilon'$-lc and $K_Y+\Delta_Y$ is numerically trivial.
Now running MMP on $-K_Y$ and using base point freeness gives another toric variety $Y'$ 
which is Fano and  $\epsilon'$-lc. By the toric version of BAB [\ref{A-L-Borisov}],  $Y'$ belongs to a bounded 
family. From this we can produce a klt strong $m$-complement 
$K_{Y'}+\Omega_{Y'}$ for some bounded $m\in\N$ which induces a klt strong $m$-complement 
$K_{Y}+\Omega_{Y}$ which in turn gives a klt strong $m$-complement 
$K_X+\Omega$. 

Finally $\Omega$ belongs to a bounded family as its coefficients are in a fixed finite set and 
its degree is bounded. This implies that $(X,\Omega+u\Lambda)$ is klt for some 
$u>0$ bounded from below. Now an easy calculation shows that the multiplicity of $T$ in 
the pullback of $\Lambda$ on $W$ is bounded from above which in turn implies the 
number of blowups in $W\to X$ is bounded as required.

\subsection{Sketch of proof of BAB}(Theorem \ref{t-BAB})
First applying [\ref{HX}, Theorem 1,3] it is enough to show that $K_X$ has a klt strong $m$-complement 
for some bounded number $m\in\N$. Running an MMP on $-K_X$ and replacing $X$ with the 
resulting model we can assume $B=0$. By Theorem \ref{t-bnd-compl-usual-local}, we know that 
we have an lc strong $n$-complement $K_X+B^+$. If $X$ is exceptional, then the complement is klt, 
so we are done in this case. To treat the general case the idea is to modify the complement 
$K_X+B^+$ into a klt one. We will do this using birational boundedness.

We need to show $\vol(-K_X)$ is bounded from above. This can be proved using arguments 
similar to the proof of the effective birationality theorem. Once we have this bound, we can 
show that $(X,B^+)$ is log birationally bounded, that is, there exist a bounded log smooth projective pair 
$(\overline{X},\Sigma_{\overline{X}})$ and  a birational map $\overline{X}\bir X$ such that 
$\Sigma_{\overline{X}}$ contains the exceptional divisors of $\overline{X}\bir X$ and the support of 
the birational transform of $B^+$. 

Next we pull back $K_X+B^+$ to a high resolution of $X$ and push it down to $\overline{X}$ 
and denote it by $K_{\overline{X}}+B^+_{\overline{X}}$. Then  $({\overline{X}},B^+_{\overline{X}})$ 
is sub-lc and $n(K_{\overline{X}}+B^+_{\overline{X}})\sim 0$. Now support of $ B^+_{\overline{X}}$ 
is contained in $\Sigma_{\overline{X}}$ so we can use the boundedness of $(\overline{X},\Sigma_{\overline{X}})$  
to perturb the coefficients of $ B^+_{\overline{X}}$. More precisely, perhaps after replacing $n$, 
there is $\Delta_{\overline{X}}\sim_\Q B^+_{\overline{X}}$ such that $({\overline{X}},\Delta_{\overline{X}})$ 
is sub-klt and $n(K_{\overline{X}}+\Delta_{\overline{X}})\sim 0$. Pulling 
$K_{\overline{X}}+\Delta_{\overline{X}}$ back to $X$ and denoting it by $K_X+\Delta$ we get a 
sub-klt $(X,\Delta)$ with $n(K_X+\Delta)\sim 0$. 

Now a serious issue here is that $\Delta$ is not necessarily effective. In fact it is by no means 
clear that its coefficients are even bounded from below. This is one of the difficult steps of the proof. 
However, this boundedness follows directly from Theorem \ref{t-bnd-lct-global}. The rest of the 
argument which modifies  $\Delta$ to get a klt complement is an easy application of complement theory.

\section{\bf Some related problems and topics} 

\subsection{Fano fibrations}\label{ss-McSh-conj} 
One of the possible outcomes of the MMP is a Mori fibre space which is an 
extremal contraction $X\to Z$ where $K_X$ is anti-ample over $Z$. This is a special kind of Fano fibration. 
Fano fibrations and more generally Fano type fibrations appear naturally in the course of applying induction 
on uniruled varieties, and in the context of moduli theory. 

Suppose now that $f\colon X\to Z$ is a Mori fibre space where $X$ is a $3$-fold with $\Q$-factorial terminal 
singularities. Mori and Prokhorov proved that if $Z$ is a surface, then $Z$ has canonical sing [\ref{Mori-Prokhorov-conic-bundle}], and 
if $Z$ is a curve, then the coefficients of the fibres of $f$ are bounded from above by $6$ [\ref{Mori-Prokhorov-del-pezzo}].

M$^{\rm c}$Kernan proposed a generalisation of the first part to higher dimension: 

\begin{conj}
Assume $d\in\N$ and $\epsilon\in \R^{>0}$. Then there is $\delta \in \R^{>0}$ such that if 
$f\colon X\to Z$ is a Mori fibre space where  $X$ is $\epsilon$-lc $\Q$-factorial of dimension $d$, 
then $Z$ is $\delta$-lc.
\end{conj}

On the other hand, independently, Shokurov proposed a more general problem which 
generalised both parts of Mori and Prokhorov result.

\begin{conj}
Assume $d\in\N$ and $\epsilon\in \R^{>0}$. Then there is $\delta \in \R^{>0}$ such that if
\begin{itemize}
\item $(X,B)$ is an $\epsilon$-lc pair of dimension $d$,  

\item $f\colon X\to Z$ is a contraction with $\dim Z>0$, 

\item $K_X+B\sim_\R 0/Z$,  and $-K_X$ is big$/Z$, 
\end{itemize}
then we can write  
$$
K_X+B\sim_\R f^*(K_Z+B_Z+M_Z)
$$
such that $(Z,B_Z+M_Z)$ is $\delta$-lc where $B_Z$ and $M_Z$ are the discriminant and moduli parts of adjunction.
\end{conj}

M$^{\rm c}$Kernan's conjecture is known in the toric case [\ref{Alexeev-Borisov}]. Shokurov's conjecture 
is known when $\dim X-\dim Z\le 1$ [\ref{B-sing-fano-fib}], in particular for surfaces, and 
open in higher dimension but we have the following general result [\ref{B-sing-fano-fib}].

\begin{thm}
Shokurov conjecture holds for those $f$ such that $(F,{{\Supp}} B|_F)$ belongs to a bounded family where 
$F$ is a general fibre of $f$.
\end{thm}

Note that by BAB (more precisely [\ref{B-BAB}, Corollary 1.2]), $F$ 
automatically belongs to a bounded family. However, one has little control over ${{\Supp}} B|_F$ 
and this is the main difficulty. This issue is similar to the difficulties which appear in the 
proof of BAB and related results. It is expected that the methods developed to prove BAB 
also works to prove Shokurov's conjecture but perhaps after some hard work.

\subsection{Minimal log discrepancies and termination}\label{ss-mld-termination}

The lc threshold plays an important role 
in birational geometry. This is clear from the proofs described in Section \ref{s-proof}. 
It is also related to the termination conjecture (\ref{conj-intro-termination}) [\ref{B-ACC-term}]. 
Another more subtle invariant of singularities is the \emph{minimal log discrepancy} (mld) also 
defined by Shorkuov. Let $(X,B)$ be a pair. The mld of $(X,B)$ denoted $\mld(X,B)$ is defined to be the minimum of log discrepancies 
$a(D,X,B)$ where $D$ runs over all prime divisors on birational model of $X$. The mld is way 
harder to treat than the lc threshold. Shokurov proposed the following:

\begin{conj}[ACC for mld's]
Assume $d\in \N$ and $\Phi\subset [0,1]$ is a set of numbers satisfying the descending chain condition (DCC). 
Then the set 
$$
\{\mld(X,B) \mid \mbox{$(X,B)$ is an lc pair and coefficients of $B$ are in $\Phi$}\}
$$
satisfies the ascending chain condition (ACC).
\end{conj}

This is known for surfaces [\ref{Alexeev-ACC}] but open in dimension $\ge 3$. Its importance is in relation with the termination 
conjecture and other topics of interest [\ref{Shokurov-mld's}][\ref{BSh}]. 
Shokurov showed that this ACC conjecture together with a semi-continuity conjecture about mld's 
due to Ambro imply the termination conjecture [\ref{Shokurov-mld's}]. The expectation is that the ACC conjecture 
can be tackled using the theory of complements and the methods described in this text but again after some hard work.

\subsection{Stable Fano varieties}\label{ss-stable-Fano}

Existence of specific metrics, e.g. K\"{a}hler-Einstein metrics, on manifolds 
is a central topic in differential geometry. Unlike canonically polarised and Calabi-Yau manifolds 
(see [\ref{Yau-78}] and references therein), Fano manifolds do not always admit such metrics. 
It is now an established fact that a Fano 
manifold admits a K\"{a}hler-Einstein metric iff it is so-called \emph{$K$-polystable} (see [\ref{CDS-stable-Fano}] 
and references therein). 

On the other hand, 
it is well-known that Fano varieties do not behave as well as canonically polarised varieties in the context 
of moduli theory. For example, the moduli space would not be separated. A 
remedy is to consider only \emph{stable} Fano's. 
The first step of constructing a moduli space is to prove a suitable boundedness result. In the smooth case this is not 
an issue [\ref{KMM-smooth-fano}] but in
the singular case boundedness is a recent result. Using methods described in Section \ref{s-proof}, 
Jiang [\ref{Jiang-stable}] proved such a result by showing that the set of $K$-semistable Fano varieties $X$ of 
fixed dimension and $\vol(-K_X)$ bounded from below forms a bounded family.

\subsection{Other topics} 

There are connections between the advances described in this text and other topics of interest  
not discussed above. Here we only mention two works very briefly.  
The papers [\ref{Lehmann-Tanimoto-Tschinkel}][\ref{Lehmann-Tanimoto-2}] 
relate boundedness of Fano's and related invariants to the geometry underlying Manin's conjecture 
on distribution of rational points on Fano varieties. 
One the other hand, [\ref{DiCerbo-Svaldi}] studies boundedness of Calabi-Yau pairs where boundedness 
of Fano varieties appear naturally.


\vspace{2cm}

\end{document}